\documentclass[12pt,oneside,reqno]{amsart}

\usepackage{hyperref}

\usepackage[headsep=30pt,footskip=30pt,twoside=false,bottom=2.5cm]{geometry}

\usepackage[all]{xy}

\usepackage{amssymb}
\usepackage{enumerate}
\usepackage{amsthm}
\usepackage{tikz-cd}
\usepackage{amsmath}
\hypersetup{
    colorlinks=true,
    linkcolor=blue,
    citecolor=red,
    filecolor=red,      
    urlcolor=violet,
}

\xymatrixcolsep{2cm}

\numberwithin{equation}{section}

\allowdisplaybreaks[1]

\usepackage{etoolbox}

\makeatletter

\patchcmd{\thesubsection}{\arabic}{\Alph}{}{}

\patchcmd{\@seccntformat}{\@secnumfont}{%
  \@secnumfont\expandafter\protect\csname format#1\endcsname}{}{}

\patchcmd{\@startsection}{\@afterindenttrue}{\@afterindentfalse}{}{}

\patchcmd{\subsection}{-.5em}{.3\linespacing}{}{}

\makeatother

\newtheorem{thm}{Theorem}[section]
\newtheorem{lem}[thm]{Lemma}
\newtheorem{prop}[thm]{Proposition}
\newtheorem{defi}[thm]{Definition}

\newtheorem{asump}[thm]{Assumption}
\newtheorem{corl}[thm]{Corollary}
\newtheorem{rem}[thm]{Remark}

\DeclareMathOperator{\Res}{{Res}}

\DeclareMathOperator{\rank}{{rank}}
\DeclareMathOperator{\Pic}{{Pic}}
\DeclareMathOperator{\Par}{{Par}}

\DeclareMathOperator{\End}{{End}}

\DeclareMathOperator{\SPar}{{SPar}}

\DeclareMathOperator{\GL}{{GL}}

\DeclareMathOperator{\pardeg}{{pardeg}}
\DeclareMathOperator{\pc}{{pc}}
\begin{document}
\baselineskip=14pt
\author{Pradeep Das}
\author{Snehajit Misra}
\author{Anoop Singh}
\address{Rajiv Gandhi Institute of Petroleum Technology, Amethi - 229304, India.}
\email[Pradeep Das]{pradeepdas0411@gmail.com}
\address{Department of Mathematics and Computing, Indian Institute of Technology (Indian School of Mines) IIT (ISM) Dhanbad - 826004, Jharkhand, India.}
\email[Snehajit Misra]{misra08@gmail.com, snehajitm@iitism.ac.in}
\address{Department of Mathematical Sciences\\ Indian Institute of Technology(BHU), Varanasi 221005, India.}
\email[Anoop Singh]{anoopsingh.mat@iitbhu.ac.in}

\title{A remark on the  Moduli of parabolic connections}

\begin{abstract}
 We consider the moduli space of stable parabolic connections of full flag type with generic rational weights on a compact Riemann surface $X$ of genus $g \ge 3$, with singularities along the parabolic points. We compute the Chow groups of the fixed determinant moduli space consisting of parabolic connections whose underlying parabolic bundle is stable.  We show that the moduli space of stable parabolic connections is not rationally connected, whereas the fixed determinant moduli space is rationally connected and in fact rational. We also prove that the determinant map to the Picard variety is universal among morphisms to abelian varieties. Finally, we show that the algebraic group of rank one holomorphic connections acts trivially on the cohomology of the moduli space of stable parabolic connections.
\end{abstract}

\maketitle

\section{Introduction and Main results}
Algebro-geometric invariants of moduli spaces play an important role in the classification of algebro-geometric objects. In this article, we consider the moduli space of stable parabolic connections with rational generic weights and full flag type. This moduli space has been constructed in  \cite{Al17,Ina13}. A compactification,  Picard group and the ring of algebraic functions of the moduli space of full-flag parabolic connections with generic rational weights have been studied in \cite{AS22}. The aim of this article is to determine the Chow groups of the moduli space of stable parabolic connections, to investigate its rationality, and to study the cohohomological aspects of the action of the algebraic group of rank one connections on it. 

Let $X$ be a compact Riemann surface with genus $g\geq 3$. Let $S = \bigl\{ x_1,x_2,\dotsc,x_m\bigr\}$ be a finite subset of $X$. Let $E$ be a holomorphic vector bundle on $X$ of rank $r$ and degree $d$. A \it quasi-parabolic structure \rm at $x\in S$ is a flag $$0=E_x^{k+1} \subsetneq E_x^{k}\subsetneq \cdots \subsetneq E_x^2\subsetneq E_x^1 = E_x$$ of linear subspaces in the fiber $E_x$. The number $k$ is the length of the flag.
Let $m_j^x =\dim_{\mathbb{C}}(E_x^j)-\dim_{\mathbb{C}}(E_x^{j+1})$ for each $j$. Then the $k$-tuple $(m_1^x,m_2^x,\dotsc,m_k^x)$ is called the type of the flag. A flag is called full if $m_i^x = 1$ for every $i$. A \it parabolic structure \rm in $E$ at $x$ is just a quasi-parabolic structure at $x$ together with a sequence of real numbers $$0\leq \alpha_1^x < \cdots <\alpha_k^x < 1,$$ called as weights. We denote by $\alpha = \bigl\{(\alpha_1^x,\dotsc,\alpha_k^x)\bigr\}_{x\in S}$ the system of real weights corresponding to a fixed parabolic structure.

A \it parabolic vector bundle \rm with a parabolic structure on $S$ is a holomorphic vector bundle $E$ together with a parabolic structure in $E$ at each point $x\in S$. We shall use the notation $E_*$ to denote a parabolic vector bundle with the underlying vector bundle $E$.

 The parabolic degree of $E_*$, denoted $\pardeg(E_*)$, is defined by
\begin{equation}
    \label{eq:para-deg}
    \pardeg(E_*) = \deg(E) + \sum\limits_{x\in S}\sum\limits_{j=1}^km_j^x\alpha_j^x,
\end{equation}
where $\deg(E)$ denotes the degree of the underlying vector bundle $E$.

We fix the integer $r\geq 1$ and $d\in \mathbb{Z}$. One can express the system of weights $\alpha = \bigl\{(\alpha_1^x,\dotsc,\alpha_r^x)\bigr\}_{x \in S}$ in matrix form $(\alpha_j^{x_i})^{1\leq i \leq m}_{1\leq j \leq r} \in \mathbb{Q}^{mr}$ for the full flag filtrations.  Inaba \cite[Definition 2.2]{Ina13} defined the system of special weights, and the system of weights which is not special is called generic. Let 
\begin{equation}
    \label{eq:W}
W^{(m)}_r(d) := \Bigl\{\alpha = (\alpha_j^{x_i})^{1\leq i \leq m}_{1\leq j \leq r} \in \mathbb{Q}^{mr} \mid \alpha  \mathrm{\,\, is \,\, generic \,\, and \,\,} d + \sum\limits_{i=1}^m\sum\limits_{j=1}^r\alpha_j^{x_i} = 0 \Bigr\}.
\end{equation}
Then the set of system of admissible parabolic weights $W^{(m)}_r(d)$ is non-empty (see \cite[Section 2.2]{BH95}).

\begin{asump}
\label{asump-1}
For fixed $r \geq 1$ and $d \in \mathbb{Z}$, we only consider those full flag parabolic weights $\alpha = \bigl\{(\alpha_1^x,\dotsc,\alpha_r^x)\bigr\}_{x \in S}$  such that $\alpha \in W^{(m)}_r(d)$, 
where the set $W^{(m)}_r(d)$ is defined in \eqref{eq:W}.
    \end{asump}

The paper is organised as follows. In Section \ref{Mod-para-conn}, we recall the definition of parabolic connections on a parabolic vector bundle. We discuss the stability and describe the moduli spaces moduli spaces $\mathcal{M}_{pc}(r,d,\alpha)$ and $\mathcal{M}_{pc}(r,d,\xi)$ of stable parabolic connections.

In Section \ref{Chow-grop}, we recall the definition of Chow groups, the localization sequence and the projective bundle formula. We use them to compute the Chow groups of the open subset $\mathcal{M}'_{pc}(r,\alpha,\xi)$  of the fixed determinant moduli space defined in \eqref{moduli-3}. 
More precise, we show the following (see Theorem \ref{thm:Chow} for the proof).

\begin{thm} Let $N = \dim(\mathcal{M}'_{pc}(r,\alpha,\xi))- \dim(\mathcal{M}(r, \alpha, \xi)) = (r^2-1)(g-1) + \frac{m}{2} (r^2-r)$. 
For every $k$ satisfying $0 \leq k \leq N$, we have a canonical isomorphism 
\begin{center}
\bf CH$_{k+(r^2-1)(g-1) + \frac{m}{2} (r^2-r)}(\mathcal{M}'_{pc}(r,\alpha,\xi)) \cong $ \bf CH$_{k} (\mathcal{M}(r, \alpha, \xi))$
\end{center}
\end{thm}

Combined with the computation of the Chow group of $1$-cycles of $\mathcal{M}(r, \alpha, \xi)$ in \cite{C21}, this gives the Chow group ${\bf CH}_{1+(r^2-1)(g-1) + \frac{m}{2} (r^2-r)}(\mathcal{M}'_{pc}(r,\alpha,\xi))$ (see Corollary~3.4) in terms of ${\bf CH}_{1} (\mathcal{M}(r, \alpha, \xi))$.

Section \ref{rational-mod} is devoted to rational connectedness and rationality of the moduli spaces $\mathcal{M}_{pc}(r,d,\alpha)$ and $\mathcal{M}_{pc}(r,d,\xi)$ of stable parabolic connections. We prove that the moduli space $\mathcal{M}_{pc}(r,d,\alpha)$ is neither rationally connected nor rational (see Proposition \ref{prop:n-rat}). And for the moduli space $\mathcal{M}_{pc}(r,\alpha,\xi)$ in \eqref{eq:mod-xi}, we prove the following  (see Theorem \ref{thm:rat} for the proof)
\begin{thm}
    The moduli space $\mathcal{M}_{pc}(r,\alpha,\xi)$ in \eqref{eq:mod-xi} is rational.
\end{thm}
There is a natural determinant map  $\det:\mathcal{M}_{pc}(r,d,\alpha)\to\operatorname{Pic}^d(X)$, we show the following  (see Theorem \ref{thm:uni}). 
\begin{thm}
    The natural morphism  $\det :  \mathcal{M}_{pc}(r,d,\alpha) \longrightarrow \Pic^d(X) $
    defined by $$(E_*,D) \mapsto \bigwedge^r E = \det(E)$$ is universal, that is, every morphism from $\mathcal{M}_{pc}(r,d,\alpha)$ to an abelian variety factors uniquely through $\det$.
\end{thm}

Finally, in  Section~\ref{rank-one-conn}, we consider the action of the algebraic group $\operatorname{Pic}^{\natural}(X)$ of rank one holomorphic connections on $\mathcal{M}_{pc}(r,d,\alpha)$ by tensor product. We prove in Theorem~\ref{thm:coh} that every element of $\operatorname{Pic}^{\natural}(X)$ acts as the identity on the cohomology $\mathrm{H}^*(\mathcal{M}_{pc}(r,d,\alpha),R)$ for any coefficient ring $R$.

\section{ Moduli space of parabolic connections}
\label{Mod-para-conn}

\subsection{Moduli of parabolic bundles}
Recall that the parabolic slope of a non-zero parabolic bundle $E_*$ on $X$ is defined as the quotient $$ \mu_{\mathrm{par}}(E_*) = \frac{\pardeg(E_*)}{\rank(E)}.$$
A subbundle and quotients of a parabolic bundle $E_*$ inherits a parabolic structure from one on $E$ in a canonical way; see \cite{MS80} for more details.

\begin{defi}
A parabolic bundle $E_*$ is said to be \it parabolic semistable \rm (resp. \it parabolic stable\rm) if for every non-zero proper parabolic subbundle $F_*$ we have 
\begin{center}
 $\mu_{\mathrm{par}}(F_*) \leq  \mu_{\mathrm{par}}(E_*)$ (resp. $<$).
\end{center}
\end{defi}

Let $\mathcal{M}^{ss}(r,d,\alpha)$ be the moduli of semi-stable parabolic vector bundles of rank $r$, deg $d$ and weight system $\alpha$ as in Assumption \ref{asump-1}. Then $\mathcal{M}^{ss}(r,d,\alpha)$ is a normal projective variety of dimension $r^2(g-1)+1+\frac{m(r^2-r)}{2}$ (see \cite{MS80}). Let $\mathcal{M}^{s}(r,d,\alpha) \subseteq \mathcal{M}^{ss}(r,d,\alpha)$ be the open subset parametrizing the stable parabolic bundles.
Since the full flag weight system is generic, we have a smooth complex projective variety
\begin{equation}
    \label{eq:moduli-1}
    \mathcal{M}(r,d,\alpha):= \mathcal{M}^{ss}(r,d,\alpha) = \mathcal{M}^s(r,d,\alpha)
\end{equation}

Let 
\begin{equation}
    \label{eq:xi}
    \xi
\end{equation}
be a fixed holomorphic line bundle over $X$ of degree $d$ such that $$d=\deg(\xi) = - \sum\limits_{x\in S}\sum\limits_{i=1}^r\alpha_i^x.$$

Let $\mathcal{M}^{ss}(r,\alpha,\xi)$ denote the moduli space of semistable parabolic vector bundles on $X$ of rank $r$ and determinant $\xi$, i.e. $\wedge^rE\cong \xi$ with weight system $\alpha$ as in Assumption \ref{asump-1}. Then $\mathcal{M}^{ss}(r,\alpha,\xi)$ is a projective variety of dimension $(r^2-1)(g-1)+\frac{m(r^2-r)}{2}$ (see \cite[p.n.~557]{BH95}). Let $\mathcal{M}^{s}(r,\alpha,\xi) \subseteq \mathcal{M}^{ss}(r,\alpha,\xi)$ be the open subset parametrizing the stable parabolic bundles. This is an irreducible smooth subvariety of $\mathcal{M}^{ss}(r,\alpha,\xi)$. Again, since the full flag weight system is generic, we have an irreducible smooth complex projective variety
\begin{equation}
    \label{eq:moduli-xi}
    \mathcal{M}(r,\alpha,\xi) := \mathcal{M}^{ss}(r,\alpha,\xi) = \mathcal{M}^s(r,\alpha,\xi).
\end{equation}

\subsection{Parabolic connections}

We now recall the notion of logarithmic connections and its residues (see \cite{BGK87}, \cite{AS22} and \cite{Del70} for details), and using these notions we will define parabolic connections. Let $E$ be a vector bundle over $X$. We denote the cotangent bundle on $X$ by $\Omega^1_X$. Let $S =x_1+x_2+\cdots+x_m$ be the reduced effective divisor associated with $S$. We define $\Omega^1_X(S):= \Omega^1_X\otimes \mathcal{O}_X(S)$.
\begin{defi}
A logarithmic connection on $E$ singular over $S$ is a $\mathbb{C}$-linear map $$D : E \longrightarrow E\otimes \Omega^1_X(S) := E \otimes \Omega^1_X \otimes \mathcal{O}_X(S)$$
which satisfies the Leibnitz identity $$D(fs) = fD(s)+df\otimes s,$$
where $f$ is a local section of $\mathcal{O}_X$ and $s$ is a local section of $E$.
\end{defi}

Given a logarithmic connection $(E, D)$ singular over $S$, we get an endomorphism denoted by 
\begin{equation}
    \label{eq:res}
    \Res(D,x_i) \in \End(E)_{x_i} = \End(E_{x_i}),
\end{equation}
 which we call the residue of the logarithmic connection $D$ at the point $x_i \in S$ (for more details see ).

\begin{defi} 
A morphism $\phi : F_* \longrightarrow E_*$ of parabolic bundles is a bundle map satisfying $$\phi(F_{x}^{i}) \subset E_{x}^{j+1}$$
whenever $\alpha(F)_{i}^x > \alpha(E)_{j}^x$ for all $x\in S$. Let $E_*$ be a parabolic vector bundle and $\phi : E \longrightarrow E$ be a morphism. We say that $\phi$ is \it strongly parabolic \it if for every $x\in S$, we have $\phi(E_x^i) \subset E_x^{i+1}$. We say that $\phi$ is \it weakly parabolic \it if it satisfies $\phi(E_x^i) \subset E_x^i$. The sheaf of strongly parabolic (resp. weakly parabolic) endomorphisms of $E_*$ is denoted by $\SPar\End(E_*) (resp. \,\,  \mathrm{Par}\End(E_*)).$ 
\end{defi}

\begin{defi}
    Let $E_*$ be a parabolic vector bundle over $X$ with a fixed system of weights $\alpha$ for the full flag filtration. A parabolic connection on $E_*$ is a logarithmic connection $D$ on the underlying vector bundle $E$ satisfying the following conditions:
\begin{enumerate}
 \item For each $x\in S$, the homomorphism induced in the filtration over the fiber $E_x$ satisfies $D\bigl(E_x^i\bigr) \subset E^i_x\otimes \Omega^1_X(S)\vert_x$ for every $i=1,2,\dotsc,r$.
 \item Since $\Res(D,x)$ preserves the filtration, it acts on each quotient $E_x^i/E_x^{i+1}$. For every $x\in S$ and for every $i=1,2,\dotsc,r$, the action of $\Res(D,x)$ on the quotient $E_x^i/E_x^{i+1}$ is the multiplication by $\alpha_i^x$. 
\end{enumerate}
We denote the parabolic connection by the pair $(E_*,D)$.  
\end{defi}

We equip $\xi$ in \eqref{eq:xi} with a parabolic connection as follows:  Let $\xi$ be given a trivial filtration at each $x\in X$ with parabolic weight $\beta^x := \sum\limits_{i=1}^r \alpha_i^x.$ Then the parabolic degree of $\xi_*$ is given by  $$\pardeg(\xi_*) = \deg(\xi)+\sum\limits_{x\in S}\beta^x = 0.$$ Now since $\pardeg(\xi_*) =0$, by \cite[Theorem 3.1]{Bis02}, $\xi_*$ admits a parabolic connection $D_{\xi_*}$ such that $\Res(D_{\xi_*},x) = \beta^x$ for every $x\in S$.  This gives a parabolic connection $(\xi_*,D_{\xi_*})$.

Given a parabolic connection $(E_*,D)$, we say $D$ is a parabolic connection for the group $\mathrm{SL}(r,\mathbb{C})$ if the logarithmic connection 
\begin{center}
 Tr$(D) : \bigwedge^rE\longrightarrow \bigwedge^rE\otimes \Omega^1_X(S)$ 
\end{center}
induced from $D$ coincide with $D_{\xi_*}$, i.e., we have an isomorphism 
$$(\bigwedge^rE_*, \mathrm{Tr}(D)) \cong (\xi_*,D_{\xi_*}).$$

\subsection{Moduli space of stable parabolic connections}
\begin{defi}
    A parabolic connection $(E_*,D)$ (for the group $\mathrm{GL}(r,\mathbb{C})$ or $ \mathrm{SL}(r,\mathbb{C})$) is said to be semi-stable (respectively, stable) if for every non-zero proper parabolic subbundle $F_*$ of $E_*$, which is invariant under $D$, i.e. $D(F) \subset F \otimes \Omega^1_X(S)$, we have 
\begin{center}
 $\mu_{\mathrm{par}}(F_*) \leq  \mu_{\mathrm{par}}(E_*)$ (respectively,  $<$ ).
\end{center}
\end{defi}
Let 
\begin{equation}
    \label{eq:moduli-2}
    \mathcal{M}_{pc}(r,d,\alpha)
\end{equation}
 be the moduli space of stable parabolic connections $(E_*,D)$ for the group $\GL(r,\mathbb{C})$ of rank $r$, degree $d$ and weight system $\alpha$ as in Assumption \ref{asump-1}.

If $(E_*, D) \in \mathcal{M}_{pc}(r,d,\alpha)$, then we have 
$$\pardeg(E_*) = \deg(E)+\sum\limits_{x\in S}\sum\limits_{i=1}^r \alpha_i^x = 0.$$
In general, the stability of $(E_*,D)$ does not imply that $E_*$ is parabolic stable.
Let 
\begin{equation}
    \label{eq:inc}
    i : \mathcal{M}'_{pc}(r,d,\alpha)\hookrightarrow \mathcal{M}_{pc}(r,d,\alpha)
\end{equation}
be the open subvariety consisting of the pairs $(E_*,D)$ whose underlying parabolic bundle $E_*$ is stable.

Let 
\begin{equation}
\label{eq:pi}
\pi : \mathcal{M}'_{pc}(r,d,\alpha) \longrightarrow \mathcal{M}(r,d,\alpha)
\end{equation}

be the forgetful map defined by sending $(E_*,D)$ to the stable parabolic bundle $E_*$.
\begin{lem}
For every $E_*\in \mathcal{M}(r,d,\alpha)$, $\pi^{-1}(E_*)$ is an affine space modelled over $H^0\bigl(X,\Omega^1_X(S)\otimes \SPar\End(E_*)\bigr)$ 
\end{lem}
\begin{proof}
    Given two parabolic connections $D$ and $D'$ on a parabolic vector bundle $E_*$, the difference $D-D'$ is an $\mathcal{O}_X$-module homomorphism from $E$ to $E \otimes \Omega^1_X$ such that the residue $$\Res(D,x)-\Res(D',x) = \Res(D-D',x)$$ acts as the zero morphism on each successive quotients $E_x^i/E_x^{i+1}$. Therefore, for every $x\in S$, we have $(D-D')(E_x^i)\subset E_x^{i+1}\otimes \Omega^1_X(S)\vert_x.$ 

Conversely, given any parabolic connection $(E_*,D)$ and 
$$\Phi\in H^0(X,\Omega^1_X(S)\otimes \SPar\End(E_*)),$$ $D+\Phi$ is again a parabolic connection on $E_*$. Thus the space $\pi^{-1}(E_*)$ of parabolic connections on $E_*$ forms an affine space modelled over the vector space $H^0(X,\Omega^1_X(S)\otimes \SPar\End(E_*)).$
\end{proof}

Let 
\begin{equation}
    \label{eq:mod-xi}
    \mathcal{M}_{pc}(r,d,\xi)
\end{equation}
be the moduli space of stable parabolic connections (for the group SL$(r,\mathbb{C})$ of rank $r$, degree $d$, generic weight system $\alpha$, and with fixed determinant $(\xi_*,D_{\xi_*}))$. The moduli space  $\mathcal{M}_{pc}(r,d,\xi)$ is a smooth irreducible quasi-projective variety of dimension $2(g-1)(r^2-1)+m(r^2-r)$. 

Let  
\begin{equation}
    \label{moduli-3}
    \mathcal{M}'_{pc}(r,d,\xi)\subseteq  \mathcal{M}_{pc}(r,d,\xi)
\end{equation}
be the Zariski open subset consisting of parabolic connections $(E_*,D)$ whose underlying parabolic vector bundle $E_*$ is parabolic stable. Let $\SPar\End'(E_*)$ and $\Par\End'(E_*)$ denote respectively the sheaves of trace zero strongly and weakly parabolic endomorphisms on $E_*$.

We have the forgetful map 
\begin{align}
\label{seq1.2}
\pi_{\xi} : \mathcal{M}'_{pc}(r,\alpha,\xi) \longrightarrow \mathcal{M}(r,\alpha,\xi)
\end{align}
which sends $(E_*,D)$ to $E_*$.  By a similar argument, one can conclude that 
for every $E_* \in \mathcal{M}(r,\alpha,\xi)$, $\pi_{\xi}^{-1}(E_*)$ is an affine space modelled over $\mathrm{H}^0\bigl(X,\Omega^1_X(S)\otimes \SPar(\End'(E_*))\bigr).$

        
\section{Chow group of the moduli space}
\label{Chow-grop}
In this section we recall the definition of Chow group of a quasi-projective scheme over a field (see \cite{F98} and \cite{V03}). Let $\mathcal{Y}$ be a quasi-projective scheme over a field $K$ of dimension $n$. We denote by $Z_k(\mathcal{Y})$ the free abelian group generated by the reduced and irreducible $k$-dimensional subvarieties of $\mathcal{Y}$. An element of $Z_k(\mathcal{Y})$ is called a $k$-dimensional algebraic cycle on $\mathcal{Y}$. Let $f\in K(\mathcal{Y})^*$. Then we have a divisor div$(f)$ on $\mathcal{Y}$ associated to the non-zero rational function $f$ on $\mathcal{Y}$. 

A $k$-cycle $\alpha\in Z_k(\mathcal{Y})$ is called rationally equivalent to zero, denoted by $\alpha \sim 0$, if there are finite number of $(k+1)$-dimensional subvarieties $W_i$ of $\mathcal{Y}$ and $f_i\in K(W_i)^*$ such that 
\begin{center}
$\alpha = \sum\limits_i $div$(f_i)$.
\end{center}
Note that the set of $k$-cycles rationally equivalent to zero form a subgroup $Z_k(\mathcal{Y})_{rat}$ of $Z_k(\mathcal{Y})$.  The quotient group 
\begin{equation}
    \label{eq:chow-gp}
    {\bf CH}_k(\mathcal{Y}) := \frac{Z_k(\mathcal{Y})}{Z_k(\mathcal{Y})_{rat}}
\end{equation}
is called the Chow group of $k$-cycles on $\mathcal{Y}$. We also consider the graded algebra 
\begin{center}
\bf CH$_*(\mathcal{Y}) = \bigoplus\limits_{k=0}^{n}$ \bf CH$_k(\mathcal{Y})$,                                          \end{center}                                        
where the graded structure is induced by the intersection products between cycles on $\mathcal{Y}$. The Chow group of $k$-cycles on $\mathcal{Y}$ with rational coefficients will be denoted by \bf CH$_k^{\mathbb{Q}}(\mathcal{Y})$\rm. We will simply write \bf CH$_k(\mathcal{Y})$ \rm when there is no confusion. We refer \cite{F98} for more details about Chow groups.

Next we mention two important lemmata which we will use to compute the Chow group of the moduli space $\mathcal{M}'_{pc}(r,\alpha,\xi)$.

Let $Y\hookrightarrow \mathcal{Y}$ be a closed subscheme, and $i :Y \longrightarrow \mathcal{Y}$ be the closed immersion, and thus $i$ is proper. Let $j : U:=\mathcal{Y}\setminus Y \longrightarrow \mathcal{Y}$ be the corresponding open immersion which is flat. Therefore we have morphisms $j^* : $ \bf CH$_k(\mathcal{Y}) \longrightarrow $ CH$_k(U)$\rm ~and $i_*:$ \bf CH$_k(Y)\longrightarrow $ CH$_k(\mathcal{Y})$ \rm of Chow groups and $j^*\circ i_* = 0$, as the cycles supported on $Y$ do not intersect $U$. Then we have the following exact sequence what is called as localization sequence :
\begin{lem} \label{lem:2.1}
 \cite[Lemma 9.12]{V03} The following sequence of abelian groups 
 \begin{center}
  \bf CH$_k(Y) \xrightarrow{i_*}$ CH$_k(\mathcal{Y}) \xrightarrow{j^*}$ CH$_k(U) \longrightarrow 0$, 
 \end{center}
is exact for every $k$ with $0\leq k \leq n = \dim(\mathcal{Y}).$
\end{lem}
\begin{lem} \label{lem:2.2}
\cite[Theorem 9.25]{V03} Let $\pi : \mathbb{P}(E)\longrightarrow \mathcal{Y}$ be the projective bundle associated to a vector bundle $E$ of rank $r$ on $\mathcal{Y}$. Then the map 
\begin{center}
 $\bigoplus\limits_{k=0}^{r-1} h^k\pi^*: \bigoplus\limits_{k=0}^{r-1}$ \bf CH$_{l-r+1+k}(\mathcal{Y})\longrightarrow$ CH$_l\bigl(\mathbb{P}(E)\bigr)$
\end{center}
is an isomorphism, where $h\in \Pic\bigl(\mathbb{P}(E)\bigr)$ denotes the class of the tautological line bundle $\mathcal{O}_{\mathbb{P}(E)}(1)$.
\end{lem}

\begin{thm} 
\label{thm:Chow}
Let $N = \dim(\mathcal{M}'_{pc}(r,\alpha,\xi))- \dim(\mathcal{M}(r, \alpha, \xi)) = (r^2-1)(g-1) + \frac{m}{2} (r^2-r)$. 
For every $k$ satisfying $0 \leq k \leq N$, we have a canonical isomorphism 
\begin{center}
\bf CH$_{k+(r^2-1)(g-1) + \frac{m}{2} (r^2-r)}(\mathcal{M}'_{pc}(r,\alpha,\xi)) \cong $ \bf CH$_{k} (\mathcal{M}(r, \alpha, \xi))$
\end{center}
\begin{proof}
Recall that the map $\pi_{\xi}$ defined in (\ref{seq1.2}) is a $T^{*}\mathcal{M}(r,\alpha,\xi)$-torsor (see \cite[Lemma 2.4]{AS22} ). Let 
$$\Psi: \mathcal{V} \longrightarrow \mathcal{M}(r,\alpha,\xi)$$ be the algebraic vector bundle such that for every Zariski open subset $U$ of $\mathcal{M}(r,\alpha,\xi)$, a section of $\mathcal{V}$ over $U$ is an algebraic function $f : \pi^{-1}_{\xi}(U) \longrightarrow \mathbb{C}$ whose restriction to each fiber $\pi^{-1}_{\xi}(E_
*)$, is an element of $\pi^{-1}_{\xi}(E_*)^{\vee}$. 

Let $(E_{*},D) \in \mathcal{M}^{\prime}_{\pc}(r,\alpha,\xi)$, and define a map $\Psi_{(E_{*},D)} : \pi^{-1}(E_{*})^{\vee} \longrightarrow \mathbb{C}$, by $$\Psi_{(E_{*},D)}(\eta) = \eta[(E_{*},D)],$$ that is in fact the evaluation map. Now, the kernel Ker$(\Psi_{(E_{*},D)})$ defines a hyperplane in $\pi^{-1}(E_{*})^{\vee}$ denoted by $H_{(E_{*},D)}$. Let $\mathbb{P}(\mathcal{V})$ be the projective bundle defined by hyperplanes in the fibre $\pi^{-1}(E_{*})^{\vee}$. Next, define a map $$i : \mathcal{M}^{\prime}_{\pc}(r,\alpha,\xi) \longrightarrow \mathbb{P}(\mathcal{V})$$ by sending $(E_{*},D)$ to the equivalence class of $H_{(E_*,D)}$, which is an open embedding.
Thus, there exists an algebraic vector bundle $$\Psi : \mathcal{V} \longrightarrow  \mathcal{M}(r,\alpha,\xi)$$ such that $\mathcal{M}'_{pc}(r,\alpha,\xi)$ is embedded in $\mathbb{P}(\mathcal{V})$ with  
$$\mathbb{H}_0 := \mathbb{P}(\mathcal{V}) \setminus \mathcal{M}'_{pc}(r,\alpha,\xi)$$ as a smooth divisor at infinity.
 Let 
 \begin{equation}
     \label{eq:proj-bun}
     \tilde{\Psi} : \mathbb{P}(\mathcal{V}) \longrightarrow \mathcal{M}(r,\alpha,\xi) 
 \end{equation}
  be the projective bundle.  Recall that for every $E_* \in \mathcal{M}(r,\alpha,\xi)$, $\pi_{\xi}^{-1}(E_*)$ is an affine space modelled over $\mathrm{H}^0\bigl(X,\Omega^1_X(S)\otimes \SPar(\End'(E_*))\bigr),$ where $\pi_\xi$ is the forgetful map defined in \eqref{seq1.2}.
  
  Let $E_*\in \mathcal{M}(r,\alpha,\xi)$, and $\theta \in \tilde{\Psi}^{-1}(E_*)\cap \mathbb{H}_0 \subset \mathbb{P}(\mathcal{V})$. Then $\theta$ represents a hyperplane in the fiber $\mathcal{V}(E_*) = \pi_{\xi}^{-1}(E_*)^{\vee}$ of the vector bundle $\mathcal{V}$. Let $H_{\theta}$ denote this hyperplane represented by $\theta$. Note that $H_{\theta} \subset \pi_{\xi}^{-1}(E_*)^{\vee}$ and $\pi_{\xi}^{-1}(E_*)$ is the affine space modelled over the vector space $H^0\bigl(X,\Omega^1_X(S)\otimes \SPar(\End'(E_*))\bigr).$ Therefore for any $f\in H_{\theta}$, we have $$df\in (T^*_{E_*} \mathcal{M}(r,\alpha,\xi))^* = T_{E_*}\mathcal{M}(r,\alpha,\xi).$$
  Since $\theta \in \mathbb{H}_0$, $H_{\theta}$ is a hyperplane, the subspace of $T_{E_*}\mathcal{M}(r,\alpha,\xi)$ generated by $\{df\}_{f\in H_{\theta}}$ is a hyperplane in $T_{E_*}\mathcal{M}(r,\alpha,\xi).$   Let $\tilde{\theta} \in \mathbb{P}(T_{E_*}\mathcal{M}(r,\alpha,\xi))$ denote this hyperplane. Thus we get a canonical isomorphism $$\mathbb{H}_0 \cong \mathbb{P}(T\mathcal{M}(r,\alpha,\xi))$$ by sending $\theta$ to $\tilde{\theta}$.
  
  Putting $Y =  \mathbb{H}_0 = \mathbb{P}(T\mathcal{M}(r,\alpha,\xi))$ and $ U= \mathbb{P}(\mathcal{V}) \setminus Y = \mathcal{M}'_{pc}(r,\alpha,\xi)$ in Lemma \ref{lem:2.1}, we get an exact sequence
  \begin{center}
   \bf CH$_l(\mathbb{P}(T\mathcal{M}(r,\alpha,\xi)))\xrightarrow{i_*}$ CH$_l(\mathbb{P}(\mathcal{V}))\xrightarrow{j^*}$ CH$_l(\mathcal{M}'_{pc}(r,\alpha,\xi))\rightarrow 0$
  \end{center}
of abelian groups for every $l = 0,1,\dotsc,\dim(\mathbb{P}(\mathcal{V})) = 2(g-1)(r^2-1)+m(r^2-r)$. 

Note that 
$$\rank(\mathcal{V}) = (g-1)(r^2-1)+\frac{m(r^2-r)}{2}+1$$
and
$$\rank(T\mathcal{M}(r,\alpha,\xi)) = \dim(\mathcal{M}(r,\alpha,\xi)) = (g-1)(r^2-1)+\frac{m(r^2-r)}{2}.$$ Hence 
$$\rank(\mathcal{V}) = \dim(\mathcal{M}(r,\alpha,\xi)) + 1.$$

For the notaional convenience, we use the following notation
\begin{equation}
    \label{moduli-xi-1}
    \mathcal{M}_\xi : = \mathcal{M}(r,\alpha,\xi)
\end{equation}

Therefore, from Lemma \ref{lem:2.2}, we get 
\begin{center}
 \bf CH$_l(\mathbb{P}(\mathcal{V})) \cong \bigoplus\limits_{i=0}^{\rank(\mathcal{V})-1}$ CH$_{l-\rank(\mathcal{V})+1+i}(\mathcal{M}_{\xi}) \cong \bigoplus\limits_{i=0}^{\dim({\mathcal{M}_{\xi})}}$ CH$_{l-\dim({\mathcal{M}_{\xi})}+2+i}(\mathcal{M}_{\xi})$
\end{center}
and 
\begin{center}
 \bf CH$_l(\mathbb{P}(T\mathcal{M}_{\xi})) \cong \bigoplus\limits_{i=0}^{\dim(\mathcal{M}_{\xi})-1}$ CH$_{l-\dim(\mathcal{M}_{\xi})+1+i}(\mathcal{M}_{\xi})$.
\end{center}
Therefore we have the following exact sequence
\begin{center}
 $\bigoplus\limits_{i=0}^{\dim(\mathcal{M}_{\xi})-1}$ \bf CH$_{l-\dim(\mathcal{M}_{\xi})+1+i}(\mathcal{M}_{\xi}) \xrightarrow{i_*} \bigoplus\limits_{i=0}^{\dim({\mathcal{M}_{\xi})}}$ CH$_{l-\dim{\mathcal{M}_{\xi}}+2+i}(\mathcal{M}_{\xi}) \xrightarrow{j^*}$ CH$_l(\mathcal{M}'_{pc}(r,\alpha,\xi))\rightarrow 0$
\end{center}
which is a short exact sequence, because $i_*$ is injective. Thus we conclude 
\begin{center}
 \bf CH$_l(\mathcal{M}'_{pc}(r,\alpha,\xi)) \cong $ CH$_{l-\dim(\mathcal{M}_{\xi})} (\mathcal{M}_{\xi})$
\end{center}
for every $l$ satisfying $\dim(\mathcal{M}_{\xi}) \leq l \leq \dim(\mathbb{P}(\mathcal{V})) = \dim(\mathcal{M}'_{pc}(r,\alpha,\xi))$.  Now, rescaling we get, for every $k$ satisfying $0 \leq k \leq  N = \dim(\mathcal{M}'_{pc}(r,\alpha,\xi))- \dim(\mathcal{M}_{\xi})$ 
\begin{center}
\bf CH$_{k+\dim(M_{\xi})}(\mathcal{M}'_{pc}(r,\alpha,\xi)) \cong $ \bf CH$_{k} (\mathcal{M}_{\xi})$
\end{center}
\end{proof}
\end{thm}

\begin{corl}
The Chow group of 1-cycles on $\mathcal{M}(r,\alpha,\xi)$ has been computed in \cite{C21}, and hence by our result we have \bf CH$_{1+\dim(\mathcal{M}(r,\alpha,\xi))}(\mathcal{M}'_{pc}(r,\alpha,\xi)) \cong $ \bf CH$_{1} (\mathcal{M}(r,\alpha,\xi))$.
\end{corl}

\section{Rationality of the moduli spaces}
\label{rational-mod}
A smooth complex variety $V$ is said to be rationally connected if any two general points on $V$ can be connected by a rational curve in $V$. Every rational variety, including the projective spaces, is rationally connected, but the converse is false.  For the theory of rational varieties, we refer to \cite{K96}.

The following
lemma is an easy consequence of the definition of rationally connectedness.
\begin{lem}
\label{lem:4.1}
Let $f : \mathcal{Y} \longrightarrow \mathcal{X}$ be a dominant rational map of complex algebraic
varieties with $\mathcal{Y}$ rationally connected. Then, $\mathcal{X}$ is rationally connected.
\end{lem}
\begin{lem}
\label{lem:4.2}
\rm(\cite[Corollary 1.3]{GHS03}) \it Let $f : \mathcal{Y} \longrightarrow \mathcal{X}$ be any dominant morphism of
complex varieties. If $\mathcal{X}$ and the general fibre of $f$ are rationally connected, then $\mathcal{Y}$ is
rationally connected.
\end{lem} 

\begin{lem} 
\label{lem:4.3}
The moduli space
$\mathcal{M}(r,d, \alpha)$ in \eqref{eq:moduli-1} is not rationally connected. 
\end{lem}
\begin{proof}
    Let $\mathcal{U}(r,d)$ be the moduli space of rank $r$ and degree $d$ stable vector bundles over $X$ such that $r$ and $d$ are coprime. Then, it is known that $\mathcal{U}(r,d)$ is a smooth projective variety.  Now, suppose that the generic weights are small (as in \cite[Definition 2.4]{AG18}), then the forgetful map $$\psi : \mathcal{M}(r,d, \alpha) \to \mathcal{U}(r,d),$$ which forgets its parabolic structure  is surjective (see \cite[Proposition 2.6]{AG18}). The fibres of $\psi$ are  flag varieties. Suppose that $\mathcal{M}(r,d, \alpha)$ is rationally connected, then in view of Lemma \ref{lem:4.1}, $\mathcal{U}(r,d)$ is rationally connected. But $\mathcal{U}(r,d)$ is not rationally connected because $\mathcal{U}(r,d)$ is birationally isomorphic to $J(X) \times \mathbb{A}^{(r^2-1)(g-1)}$ (\cite[Theorem 1.1]{KS99}) which is a contradiction. Thus, $\mathcal{M}(r,d, \alpha)$ is not rationally connected for small weights $\alpha$. Next, we show that it is true for any generic weights.  For instance,   the two moduli space $\mathcal{M}(r,d, \beta)$ and $\mathcal{M}(r,d, \gamma)$ are birational for any two generic weights $\beta, \gamma \in W^{(m)}_r(d)$ (follows from \cite[Theorem 4.2 and Proposition 4.3]{BY99}). Since rationally connectedness is a birational property.  Therefore, $\mathcal{M}(r,d, \alpha)$ is not rationally connected.
\end{proof}

\begin{prop} 
\label{prop:n-rat}
The moduli space 
$\mathcal{M}_{pc}(r,d,\alpha)$ in \eqref{eq:moduli-2} is not rationally connected, and hence not rational.
\end{prop}
\begin{proof}
    We will show that $\mathcal{M}'_{pc}(r,d,\alpha)$ in \eqref{moduli-3} is not rationally connected, which is an open dense subset of $\mathcal{M}_{pc}(r,d,\alpha)$. We consider the forgetful map $\pi$ defined in \eqref{eq:pi}. Assume that $\mathcal{M}'_{pc}(r,d,\alpha)$ is rationally connected. The fibre of $\pi$ in \eqref{eq:pi} is an affine space and hence rational and hence rationally connected. Therefore, $\pi$ is a dominant rational map. Then, in view of  Lemma \ref{lem:4.1}, $\mathcal{M}(r,d, \alpha)$ is rationally connected, which contradicts Lemma \ref{lem:4.3}. 
\end{proof}

\begin{prop}
\label{prop:RC-1}
The moduli space $\mathcal{M}_{pc}(r,\alpha,\xi)$ in \eqref{eq:mod-xi} is rationally connected.
\begin{proof}
Recall that we have the forgetful map (in \ref{seq1.2}) $$\pi_{\xi} : \mathcal{M}'_{pc}(r,\alpha,\xi) \longrightarrow \mathcal{M}(r,\alpha,\xi)$$ which sends $(E_*,D)$ to $E_*$.  By a similar argument, one can conclude that 
for every $E_* \in \mathcal{M}(r,\alpha,\xi)$, $\pi_{\xi}^{-1}(E_*)$ is an affine space modelled over $H^0\bigl(X,\Omega^1_X(S)\otimes \SPar(\End'(E_*))\bigr).$ 
Since $ \mathcal{M}(r,\alpha,\xi)$ is rationally connected by using \cite{BY99}, from Lemma \ref{lem:4.2}, we have $\mathcal{M}'_{pc}(r,\alpha,\xi)$ is rationally connected. Now the proposition follows from the fact that rationally connectedness is a birational invariant,
and $\mathcal{M}'_{pc}(r,\alpha,\xi)$  is a dense open subset of $\mathcal{M}_{pc}(r,\alpha,\xi)$.
 \end{proof}
\end{prop}

In fact, we have a stronger result as follows proved by using a different technique.
\begin{thm}
    \label{thm:rat}
    The moduli space $\mathcal{M}_{pc}(r,\alpha,\xi)$ in \eqref{eq:mod-xi} is rational.
\end{thm}

\begin{proof}
    From \cite[Theorem 6.2]{BH95} it follows that the moduli space  $\mathcal{M}(r,\alpha,\xi)$ is rational. 
    Consider the projective bundle $$\tilde{\Psi} : \mathbb{P}(\mathcal{V}) \longrightarrow \mathcal{M}(r,\alpha,\xi) $$
    constructed in \eqref{eq:proj-bun} associated to the vector bundle $\mathcal{V}$ over $\mathcal{M}(r,\alpha,\xi)$.  Since the vector bundle is Zariski locally trivial, there is a nonempty Zariski open subset
$U \subset \mathcal{M}(r,\alpha,\xi)$ such that
$$
{\tilde{\Psi}}^{-1}(U)\ \cong\ U \times\mathbb P^m, 
$$
where $m = \operatorname{rk}\mathcal{V} - 1$. Let $k = \mathrm{dim}_{\mathbb{C}} \mathcal{M}(r,\alpha,\xi)$, there is a birational map   
$
U\dashrightarrow\mathbb{P}^k.
$
Since $\mathbb{P}(\mathcal{V})$ is irreducible and
${\tilde{\Psi}}^{-1}(U)$ is a dense open subset, we have the following composition of birational maps
$$
\mathbb{P}(\mathcal{V})
\dashrightarrow{\tilde{\Psi}}^{-1}(U)\ \cong\ U \times\mathbb P^m
\dashrightarrow \mathbb{P}^k
\times\mathbb{P}^m.
$$
Since $\mathbb{P}^k\times\mathbb{P}^m$
is rational, $\mathbb{P}(\mathcal{V})$ is rational. Therefore,
its nonempty open subset  $\mathcal{M}'_{pc}(r,\alpha,\xi) \subset \mathbb{P}(\mathcal{V})$
is rational. Since the open subset $\mathcal{M}'_{pc}(r,\alpha,\xi)$ is dense in $\mathcal{M}_{pc}(r,\alpha,\xi)$,
therefore, the moduli space $\mathcal{M}_{pc}(r,\alpha,\xi)$ is rational. This completes the proof.
\end{proof}

Recall that  $\mathcal{M}_{pc}(r,d,\alpha)$ is the moduli space parameterizing stable parabolic connections $(E_*,D)$ for the group $\GL(r,\mathbb{C})$ of rank $r$, degree $d$ and weight system $\alpha$ such that $$\pardeg(E_*) = \deg(E)+\sum\limits_{x\in S}\sum\limits_{i=1}^r \alpha_i^x = d+ \sum\limits_{x\in S}\sum\limits_{i=1}^m \alpha_i^x = 0.$$ 

\begin{defi}
Let $\mathcal{X}$ be a complex algebraic variety and $A_0$  a  complex abelian variety. A surjective morphism $$p_0: \mathcal{X} \longrightarrow A_0,$$
 is called  universal  if for any algebraic morphism $p : \mathcal{X} \longrightarrow A,$ where $A$ is any complex abelian variety, there is a unique morphism
$$\gamma : A_0 \longrightarrow A$$
such that $\gamma\circ p_0 = p$.
\end{defi}

\begin{thm}
\label{thm:uni}
    The natural morphism 
    \begin{equation}
        \label{eq:det}
        \det :  \mathcal{M}_{pc}(r,d,\alpha) \longrightarrow \Pic^d(X) 
    \end{equation}
    defined by $$(E_*,D) \mapsto \wedge^r E = \det(E)$$ is universal.
\end{thm}
    
    \begin{proof}
First we show that the determinant map in \eqref{eq:det} is surjective.  Let $z\in\operatorname{Pic}^{d}(X)$. 
Since $z\otimes \xi^{-1}\in\operatorname{Pic}^{0}(X)$, where $\xi$ is in \eqref{eq:xi}
and the multiplication by $r$-morphism ($\eta \mapsto \eta^{\otimes r}$) $$[r]\colon\operatorname{Pic}^{0}(X)\longrightarrow
\operatorname{Pic}^{0}(X)$$
is surjective, there exists a holomorphic line bundle $L_0$ of degree zero
such that $L_0^{\otimes r}\cong \xi\otimes z^{-1},$ which gives  $z\otimes L_0^{\otimes r}\cong \xi.$ Now, 
consider the moduli space $\mathcal{M}_{pc}(r,d,\xi)$ in \eqref{eq:mod-xi} and choose a point $(E_*,D)\in 
\mathcal{M}_{pc}(r,d,\xi)$. Since $\deg L_0=0$, the line bundle $L_0$ admits a holomorphic connection and 
 after fixing the above isomorphism $z\otimes L_0^{\otimes r}\cong\xi,$ there is a unique holomorphic connection 
 $\nabla$ on $L_0$ such that the connection $D\otimes\operatorname{Id}_{L_0}+\operatorname{Id}_{E}\otimes\nabla$
on $E\otimes L_0$ induces the prescribed determinant connection
$D_{\xi_*}$ on $\xi_*$. Indeed, the two connections on the determinant line have the same residues
$\beta_x$ at every $x\in S$, and hence their difference is a global
holomorphic $1$-form. Dividing this difference by $r$ uniquely determines
$\nabla$. The parabolic structure on $E_*\otimes L_0$ has the same full flag and the
same weights $\alpha$, because $L_0$ has degree zero and is equipped with the
trivial parabolic structure.
Since tensoring with a degree-zero line bundle $L_0$ preserves parabolic
stability, we obtain
\[
\left(
E_*\otimes L_0,\,
D\otimes\operatorname{Id}_{L_0}
+
\operatorname{Id}_{E}\otimes\nabla
\right)
\in \mathcal{M}_{pc}(r,d,\alpha)
\]
with 
its determinant 
\[
\det(E\otimes L_0)
=
\det(E)\otimes L_0^{\otimes r}
\cong
\xi\otimes(\xi\otimes z^{-1})^{-1}
\cong z.
\]

Next we show that the fibers of the determinant morphism $\det$ in \eqref{eq:det} are rationally connected. 
Fix $z\in\operatorname{Pic}^{d}(X)$
and set
\begin{equation}
    \label{eq:fiber}
    F_z:={\det}^{-1}(z).
\end{equation}
Let $z_*$ denote the parabolic line bundle obtained by giving $z$ the
weight $\beta_x$ at each $x\in S$. Since 
$$ \deg(z)+\sum_{x\in S}\beta_x = d+\sum_{x\in S}\sum_{i=1}^r\alpha_i^x = 0, $$
by \cite[Theorem 3.1]{Bis02}, $z_*$ admits a parabolic connection $D_{z_*}$ such that $\Res(D_{z_*},x) = \beta^x$ for every $x\in S$. 
Let
\begin{equation}
    \label{eq:mod-z}
    \mathcal{M}_z :=
\mathcal{M}_{pc}(r,d,z)
\end{equation}
denote the moduli space of stable parabolic connections with fixed
determinant connection $(z_*,D_{z_*})$. We first note that $\mathcal{M}_z$ in \eqref{eq:mod-z} is isomorphic to
$\mathcal{M}_{pc}(r,d,\xi)$ in \eqref{eq:mod-xi}. Indeed, choose $L_0\in\operatorname{Pic}^{0}(X)$ with
$L_0^{\otimes r}\cong \xi\otimes z^{-1}$ and as above choose the unique holomorphic connection $\nabla_0$ on $L_0$
for which tensorization by $(L_0,\nabla_0)$ identifies the determinant
connection $(z_*,D_{z_*})$ with $(\xi_*,D_{\xi_*})$. This gives an algebraic
isomorphism
\begin{equation}
    \label{eq:iso-z}
    \mathcal{M}_z \stackrel{\sim}{\longrightarrow}
\mathcal{M}_{\mathrm{pc}}(r,d,\xi).
\end{equation}
Consequently, in view of Proposition \ref{prop:RC-1},  the moduli space $\mathcal{M}_z $ in \eqref{eq:mod-z} is rationally connected. We next describe $F_z$. Given $(E_*,D)\in F_z,$
the induced connection on $\det(E) \cong z$ has residue $\beta_x$ at each
$x\in S$. Hence $\operatorname{Tr}(D)-D_{z_*}$
has no logarithmic poles and therefore determines an element 
\begin{equation}
    \label{eq:omega}
    \omega_{(E_*,D)}
\in \mathrm{H}^0(X,\Omega_X^1).
\end{equation}
Define $D^0
:=
D-\frac{1}{r}\omega_{(E_*,D)}
\otimes\operatorname{Id}_E.$
Then
\[
\operatorname{Tr}(D^0)
=
\operatorname{Tr}(D)-\omega_{(E_*,D)}
=
D_{z_*}.
\]
Thus $(E_*,D^0)\in \mathcal{M}_z$ and therefore we obtain a map
\begin{equation}
    \label{eq:pi-z}
    \Pi_z \colon F_z \longrightarrow \mathcal{M}_z \times \mathrm{H}^0(X,\Omega_X^1)
\end{equation}
defined by $\Pi_z (E_*, D) = \left((E_*,D^0),\omega_{(E_*,D)}\right).$ The homomorphism $\Pi_z$ in \eqref{eq:pi-z}
is an algebraic isomorphism. Indeed, given $ \left( (E_*, D), \omega \right) \in \mathcal{M}_z \times \mathrm{H}^0(X,\Omega_X^1)$, we define $D_\omega = D + \frac{1}{r} \omega \otimes \operatorname{Id}_E.$
Then $D_\omega$ is again a parabolic connection $E_*$ with the same parabolic
weights, and its induced determinant connection is $\operatorname{Tr}(D_\omega)
=D_z+\omega.$ In particular, the underlying determinant line bundle is still $z$, so $(E_*,D_\omega)\in F_z.$
Since $\mathcal{M}_z$ is rationally connected and $\mathrm{H}^0 (X, \Omega^1_X)$ is an affine space, $F_z$ is rationally connected.

To prove that $\det$ in \eqref{eq:det} is universal, let $A$ be a complex abelian variety and $$p : \mathcal{M}_{pc}(r,d,\alpha) \longrightarrow A $$ be an algebraic morphism. fix a point $z\in \Pic^d(X)$ and consider the restriction map 
$$ p_z := p_{\vert {\det}^{-1}(z)} = p_{\vert F_z} : F_z \longrightarrow A.$$ Recall that every algebraic morphism from a rationally connected variety to an
abelian variety is constant. As observed above $F_z$ is rationally connected, so $p_z$ is constant. 
Consequently, we get a set theoretic map  $ \gamma : \operatorname{Pic}^{d}(X)\longrightarrow A$
such that 
\begin{equation}
    \label{eq:gamma}
    p=\gamma\circ\det.
\end{equation}
Now, we show that $\gamma$ is an algebraic morphism. Since the moduli space $\mathcal{M}_{pc}(r,d,\alpha)$ is smooth and $\operatorname{Pic}^{d}(X)$ is smooth, it follows that $\det$ in \eqref{eq:det} is smooth and hence faithfully flat. Consider now the fibre product
$ M := \mathcal{M}_{pc}(r,d,\alpha) \times_{\operatorname{Pic}^{d}(X)} \mathcal{M}_{pc}(r,d,\alpha)$ equipped with two projections $q_1$ and $q_2$. For every point $(u,v) \in M$ of this fibre product, we have $\det(u)=\det(v),$ which implies $\gamma (\det (u)) = \gamma (\det (v))$ and hence using \eqref{eq:gamma}, we get $p(u) = p(v)$. Consequently, $p\circ q_1=p\circ q_2.$
Since $\det$ is faithfully flat and $A$ is separated, the fpqc descent
criterion for morphisms gives a unique algebraic morphism $\gamma' \colon
\operatorname{Pic}^{d}(X)\longrightarrow A$
such that
\[
p=\gamma'\circ\det.
\]
From uniqueness of $\gamma'$, we get that $\gamma = \gamma'.$ This completes the proof.
\end{proof}

\section{Action of rank one connection on the moduli space}
\label{rank-one-conn}

Let $\Pic^\natural(X)$ be the moduli space of pairs $(N,\nabla)$, with $N$ a
holomorphic line bundle on $X$ and $\nabla$ a holomorphic connection on $N$.
Since $\dim X=1$, such $\nabla$ is automatically flat, and $\deg N=0$.
It is a commutative complex algebraic group under $\otimes$, with unit
$(\mathcal O_X,d)$. By Atiyah-Weil criterion, a holomorphic line bundle $N$ over a compact connected Riemann surface $X$ of degree zero admits a holomorphic connection. 

By the Riemann--Hilbert correspondence, $\Pic^\natural(X)$ is analytically
isomorphic to $$\operatorname{Hom}(\pi_1(X),\mathbb C^*)\cong(\mathbb C^*)^{2g}, $$
which shows that $\Pic^\natural(X)$ is path-connected in analytic topology. 

Note that if $(E_*, D)$ is a stable parabolic connections of type $(r, d, \alpha)$, then for any $(N, \nabla) \in \Pic^\natural(X)$, the connection $(E_*\otimes N, D\otimes\operatorname{Id}_{N}+\operatorname{Id}_{E}\otimes\nabla)$ is a stable
parabolic connection of type $(r,d,\alpha)$.

Thus, for $(N,\nabla)\in\Pic^\natural(X)$, we get a morphism
\begin{equation}
    \label{eq:mor-phi}
\Phi_{(N, \nabla)} : \mathcal{M}_{pc}(r,d,\alpha) \longrightarrow \mathcal{M}_{pc}(r,d,\alpha)
\end{equation}
defined by 
$$\Phi_{(N,\nabla)}(E_*,D) = (E_*\otimes N, D\otimes\operatorname{Id}_{N}+\operatorname{Id}_{E}\otimes\nabla).$$
Note that $\Phi_{(N,\nabla)}$ is an automorphism of $\mathcal{M}_{pc}(r,d,\alpha)$ with the inverse $\Phi_{(N^{-1}, \nabla^{-1})}$, where $\nabla^{-1}$ is the dual connection on $N^{-1}$.

In fact, we get an action of $\Pic^\natural(X)$ on the moduli space $\mathcal{M}_{pc}(r,d,\alpha)$ described as follows: define 
\begin{equation}
    \label{eq:action-mu}
    \mu : \Pic^\natural(X) \times \mathcal{M}_{pc}(r,d,\alpha) \longrightarrow \mathcal{M}_{pc}(r,d,\alpha)
\end{equation}
by
$$\mu ( (N, \nabla), (E_*, D)) = \Phi_{(N,\nabla)}(E_*,D).$$
It is easy to see that $\mu$ in \eqref{eq:action-mu} and hence $\Phi_{(N, \nabla)}$ in \eqref{eq:mor-phi} are algebraic morphisms between complex  algebraic varieties.

For any commutative ring $R$, the automoprhism $\Phi_{(N, \nabla)}$ in \eqref{eq:mor-phi} induces a map on singular cohomology 
\begin{equation}
    \label{eq:sing-coh}
    \Phi_{(N, \nabla)}^* : \mathrm{H}^*(\mathcal{M}_{\mathrm{pc}}(r,d,\alpha), R) \longrightarrow \mathrm{H}^*(\mathcal{M}_{\mathrm{pc}}(r,d,\alpha),R).
\end{equation}

\begin{thm}
\label{thm:coh}
    The induced morphism $\Phi_{(N, \nabla)}^*$ in \eqref{eq:sing-coh} on the cohomology  is the identity map. 
\end{thm}

\begin{proof}
    We show that the automorphism $\Phi_{(N,\nabla)}$ in \eqref{eq:mor-phi} is homotopic to the identity automorphism $ \operatorname{Id}_0 := \operatorname{Id}_{\mathcal{M}_{\mathrm{pc}}(r,d,\alpha)}$. Since 
    $\Pic^\natural(X)$ is path-connected in analytic topology, there exists a continuous path $\sigma : [0, 1] \longrightarrow \Pic^\natural(X)$ such that $\sigma (0) = (\mathcal O_X,d)$ and $\sigma (1) = (N, \nabla)$.
    Define $$ H : [0, 1] \times \mathcal{M}_{\mathrm{pc}}(r,d,\alpha) \longrightarrow \mathcal{M}_{\mathrm{pc}}(r,d,\alpha)$$ by $H(t, y) = \mu (\sigma(t), y)$, where $\mu$ is the algebraic action defined in \eqref{eq:action-mu}. Since $\mu$ in \eqref{eq:action-mu} is a morphism of varieties, hence continuous, and $\sigma$ is a continuous path, so $H$ is continuous.
We have $H(0,\cdot) = \operatorname{Id}_0 $ and
$H(1,\cdot) = \Phi_{(N,\nabla)}$. Since, homotopic maps induce the same map in singular
cohomology with any coefficients, so $\Phi_{(N,\nabla)}^* = \operatorname{Id}^*_0 = \operatorname{Id}_{\mathrm{H}^*(\mathcal{M}_{\mathrm{pc}}(r,d,\alpha), R)}$.
\end{proof}

\begin{rem} If we consider the automorphism of the fixed determinant moduli space $\mathcal{M}_{\mathrm{pc}}(r,d,\xi)$,  $\Phi_{(N,\nabla)}$
preserves the determinant locus only when $(N,\nabla)^{\otimes r}\cong(\mathcal O_X,d)$,
that is for $(N,\nabla)$ in the finite subgroup $\Pic^\natural(X)[r]$. The
connectedness argument therefore does not apply, and hence the above Theorem \ref{thm:coh} need not be true in general for the moduli space $\mathcal{M}_{\mathrm{pc}}(r,d,\xi)$.
\end{rem}


\section*{acknowledgement}
   A.S. is partially supported by ANRF/ARGM/2025/000670/MTR.

\end{document}